\documentclass{amsart}
\usepackage{amscd,amssymb,amsxtra}

\theoremstyle{plain}
\newtheorem{Thm}{Theorem}[section]
\newtheorem{Prop}[Thm]{Proposition}
\newtheorem{Lem}[Thm]{Lemma}
\newtheorem{Cor}[Thm]{Corollary}
\newtheorem{Claim}[Thm]{Claim}

\theoremstyle{definition}
\newtheorem{Def}[Thm]{Definition}
\newtheorem{Question}[Thm]{Question}
\newtheorem{Conj}[Thm]{Conjecture}
\newtheorem{Hyp}[Thm]{Hypothesis}

\theoremstyle{remark}

\newtheorem{Remark}[Thm]{Remark}

\numberwithin{equation}{section}

\newcommand{\dom}{\operatorname{dom}}
\newcommand{\res}{\upharpoonright}

\newcommand{\Aut}{\operatorname{Aut}}
\newcommand{\Inn}{\operatorname{Inn}}
\newcommand{\pred}{\operatorname{pred}}
\newcommand{\hit}{\operatorname{ht}}
\newcommand{\cf}{\operatorname{cf}}
\newcommand{\ran}{\operatorname{ran}}
\newcommand{\Lev}{\operatorname{Lev}}
\newcommand{\Sym}{\operatorname{Sym}}
\newcommand{\Wr}{\text{wr}}

\begin{document}

\title{Changing the heights of automorphism towers}
\author{Joel David Hamkins}
\address
{Mathematics Department \\
College of Staten Island \\
City University of New York \\
Staten Island \\
New York 10314}
\author{Simon Thomas}
\address
{Mathematics Department \\
Rutgers University \\
New Brunswick \\
New Jersey 08903}
\thanks{The research of the second author was partially supported
by NSF Grants.}


\begin{abstract}
If $G$ is a centreless group, then $\tau(G)$ denotes the height
of the automorphism tower of $G$. We prove that it is consistent that for
every cardinal $\lambda$ and every ordinal $\alpha < \lambda$, there exists
a centreless group $G$ such that \\
(a) $\tau(G) = \alpha$; and \\
(b) if $\beta$ is any ordinal such that
$1 \leq \beta < \lambda$, then there exists a notion of
forcing $\mathbb{P}$, which preserves cofinalities and cardinalities,
such that $\tau(G) = \beta$ in the corresponding generic extension
$V^{\mathbb{P}}$.
\end{abstract}

\maketitle
\section{Introduction} \label{S:intro}
If $G$ is a centreless group,
then there is a natural embedding $e_{G}$ of $G$ into its automorphism
group $\Aut (G)$, obtained by sending each $g\in G$ to the
corresponding inner automorphism $i_{g} \in\Aut (G)$. In this
paper, we shall always work with the left action of $\Aut (G)$
on $G$.  Thus $i_{g}(x) = g x  g^{-1}$ for all $x \in G$.  
If $\pi \in \Aut (G)$ and $g \in
G$, then $\pi i_{g} \pi^{-1} = i_{\pi(g)}$.  Hence the group of inner
automorphisms $\Inn (G)$ is a normal subgroup of $\Aut (G)$.  Also
$C_{\Aut (G)} (\Inn (G)) = 1$.  In particular, $\Aut (G)$ is a centreless group.  
This enables
us to define the automorphism tower of $G$ inductively as follows.

\begin{Def} \label{D:tower}
\begin{enumerate}
\item[(a)]  $G_0 = G$.
\item[(b)]  Suppose that $G_{\alpha}$ has been defined. Then $G_{\alpha + 1}$
is chosen to be a group such that
\begin{enumerate}
\item[(i)]  $G_{\alpha} \leqslant G_{\alpha+1}$; and
\item[(ii)]  there exists an isomorphism $\varphi$ such that the
following diagram commutes.
\[
\begin{CD}
  G_{\alpha+1}  @>\varphi>> \Aut (G_{\alpha})\\
    @A\text{inc}AA                           @A \text{inc}AA\\
G_{\alpha} @> e_{G_{\alpha}}>>  \Inn (G_\alpha)
\end{CD}
\]
\end{enumerate}
(Here {\em inc} denotes the inclusion map.  This corresponds to identifying 
$G_{\alpha}$ with $\Inn (G_{\alpha})$. There is actually a unique
such isomorphism $\varphi$.  This allows us to speak of {\em the}
automorphism tower of $G$.)
\item[(c)]  If $\lambda$ is a limit ordinal, then $G_{\lambda} = 
\underset{ \alpha < \lambda}{\bigcup} G_{\alpha}$.
\end{enumerate}
\end{Def}

The automorphism tower is said to terminate if there exists an ordinal
$\alpha$ such that $G_{\beta} = G_{\alpha}$ for all $\beta > \alpha$.
This occurs if and only if there exists an ordinal $\alpha$ such that
$G_{\alpha}$ is a complete group. (A centreless group G is said to
be {\em complete\/} if $\Aut (G) = \Inn (G)$.) A
classical result of Wielandt \cite{w} says that if $G$ is finite, then the
automorphism tower terminates after finitely many steps.  In 
\cite{t1}, it was shown that the automorphism tower of an arbitrary
centreless group eventually terminates; and that for each ordinal $\alpha$,
there exists a group whose automorphism tower terminates in exactly
$\alpha$ steps.

\begin{Def} \label{D:tau}
If $G$ is a centreless group, then the {\em height\/} $\tau(G)$
of the automorphism tower of $G$ 
is the least ordinal $\alpha$ such that $G_{\beta} = G_{\alpha}$
for all $\beta > \alpha$.
\end{Def}

Let $V$ denote the ground model, and let $G \in V$ be a centreless
group. If $M$ is a generic extension of $V$, then $\tau^{M}(G)$ denotes
the value of $\tau(G)$, when $\tau(G)$ is computed within $M$. In \cite{t2},
it was shown that there exist a centreless group $G \in V$ and a $c.c.c.$
notion of forcing $\mathbb{P}$ such that $\tau(G) = 0$ and
$\tau^{V^{\mathbb{P}}}(G) \geq 1$. This is not a very surprising result.
It was to be expected that there should be a complete group $G$ which
possessed an outer automorphism in some generic extension $V^{\mathbb{P}}$.
More surprisingly, it was also shown in \cite{t2} that there exists a
centreless group $G$ such that
\begin{enumerate}
\item[(a)] $\tau(G) = 2$; and
\item[(b)] if $\mathbb{P}$ is any notion of forcing which adjoins
a new real, then $\tau^{V^{\mathbb{P}}}(G) = 1$.
\end{enumerate}
Thus the height of the automorphism tower of a centreless group $G$
may either increase or decrease in a generic extension. In fact, if
$M$ is a generic extension, then it is difficult to think of any
constraints on $\tau(G)$ and $\tau^{M}(G)$; apart from the obvious
one that if $\tau(G) \geq 1$, then $\tau^{M}(G) \geq 1$. (If $G$
possesses an outer automorphism $\pi \in \Aut (G) \smallsetminus
\Inn (G)$ in $V$, then $\pi$ remains an outer automorphism in $M$.)
These considerations led the second author to make the following
conjecture in \cite{t2}.

\begin{Conj} \label{C:first}
Let $\alpha$, $\beta$ be ordinals such that if $\alpha \geq 1$,
then $\beta \geq 1$. Then there exist a centreless group $G$ 
and a notion of forcing $\mathbb{P}$ such that $\tau(G) = \alpha$
and $\tau^{V^{\mathbb{P}}}(G) = \beta$.
\end{Conj}

In this paper, we will prove the consistency of a substantial
strengthening of Conjecture \ref{C:first}.

\begin{Thm} \label{T:main}
It is consistent that for
every infinite cardinal $\lambda$ and 
every ordinal $\alpha < \lambda$, there exists
a centreless group $G$ with the following properties. 
\begin{enumerate}
\item[(a)] $\tau(G) = \alpha$. 
\item[(b)] If $\beta$ is any ordinal such that
$1 \leq \beta < \lambda$, then there exists a notion of
forcing $\mathbb{P}$, which preserves cofinalities and cardinalities,
such that $\tau^{V^{\mathbb{P}}}(G) = \beta$.
\end{enumerate}
\end{Thm}

It should be pointed out that this is {\em not\/} the strongest
conceivable consistency result on the nonabsoluteness of the heights
of automorphism towers. By \cite{t2}, if $G$ is an infinite centreless
group, then the automorphism tower of $G$ terminates in strictly
less than $\left( 2^{|G|} \right)^{+}$ steps. However,
$\left( 2^{|G|} \right)^{+}$ can be an arbitrarily large cardinal
in generic extensions of the ground model $V$. Thus the following
problem remains open.

\begin{Question} \label{Q:noway}
Does there exist a complete group $G$ such that for every ordinal
$\alpha$, there exists a notion of forcing $\mathbb{P}$, which
preserves cofinalities and cardinalities, such that
$\tau^{V^{\mathbb{P}}}(G) = \alpha$?
\end{Question}

In Section \ref{S:norm}, we will present an essentially algebraic
argument which shows that Theorem \ref{T:main} is a consequence
of the following result.

\begin{Thm} \label{T:graph}
It is consistent that for every regular cardinal $\kappa \geq \omega$, there
exists a set $\{ \Gamma_{\alpha} \mid \alpha < \kappa^{+} \}$ of 
pairwise nonisomorphic connected rigid graphs with the
following property. If  $E$ is any equivalence relation on
$\kappa^{+}$, then there exists a notion of forcing
$\mathbb{P}$ such that
\begin{enumerate}
\item[(a)] $\mathbb{P}$ preserves cofinalities and cardinalities;
\item[(b)] $\mathbb{P}$ does not adjoin any new $\kappa$-sequences 
of ordinals;
\item[(c)] each graph $\Gamma_{\alpha}$ remains rigid in $V^{\mathbb{P}}$;
\item[(d)] $\Gamma_{\alpha} \simeq \Gamma_{\beta}$ in $V^{\mathbb{P}}$ if{f}
$\alpha \mathrel{E} \beta$.
\end{enumerate}
\end{Thm}

Here a structure $\mathcal{M}$ is said to be {\em rigid\/} if
$\Aut (\mathcal{M}) = \{ id_{\mathcal{M}} \}$. Thus clauses
(\ref{T:graph})(c) and (\ref{T:graph})(d) imply that if 
$\alpha \mathrel{E} \beta$, then there exists
a {\em unique\/} isomorphism 
$\pi : \Gamma_{\alpha} \to \Gamma_{\beta}$ in $V^{\mathbb{P}}$.
Theorem \ref{T:graph} will be proved in Section \ref{S:rigid}.

Our set-theoretic notation mainly follows that 
of Jech \cite{j2}. Thus if $\mathbb{P}$
is a notion of forcing and $p$,$q \in \mathbb{P}$, then $q \leq p$
means that $q$ is a strengthening of $p$. We say that $\mathbb{P}$
is $\kappa$-closed if for every $\lambda \leq \kappa$, every 
descending sequence of elements of $\mathbb{P}$
\[
p_{0} \geq p_{1} \geq \dots \geq p_{\xi} \geq \dots , \qquad  \xi < \lambda ,
\]
has a lower bound in $\mathbb{P}$. If $V$ is the ground model,
then we will denote the generic extension by $V^{\mathbb{P}}$ if
we do not wish to specify a particular generic filter 
$H \subseteq \mathbb{P}$. If we want to emphasize that the term $t$
is to be interpreted in the generic extension $M$, then we write
$t^{M}$. For example, if $\Gamma \in V$ is a graph, then
$\Aut^{M}(\Gamma)$ denotes the automorphism group of
$\Gamma$, when the automorphism group is computed in $M$. 
The class of all ordinals will be denoted by $On$.

Our group-theoretic notation is standard. For example, if $G$ is a group,
then $Z(G)$ denotes the centre of $G$. A {\em permutation group\/}
is a pair $\left( G, \Omega \right)$, where $G$ is a subgroup of
$\Sym (\Omega)$. A pair $(f, \varphi )$ is a 
{\em permutation group isomorphism\/} from $\left( G, \Omega \right)$ onto
$\left( H, \Delta \right)$ if the following conditions are satisfied.
\begin{enumerate}
\item[(i)] $f : G \to H$ is a group isomorphism.
\item[(ii)] $\varphi : \Omega \to \Delta$ is a bijection.
\item[(iii)] For all $g \in G$ and $x \in \Omega$,
$f(g)(\varphi(x)) = \varphi( g(x))$.
\end{enumerate}

\section{Normaliser towers} \label{S:norm}
In this section, we will show that Theorem \ref{T:main} is a consequence
of Theorem \ref{T:graph}. So throughout this section, we will assume
that the following hypothesis holds in the ground model $V$.

\begin{Hyp} \label{H:graph}
For every regular cardinal $\kappa \geq \omega$, there
exists a set $\{ \Gamma_{\alpha} \mid \alpha < \kappa^{+} \}$ of 
pairwise nonisomorphic connected rigid graphs with the
following property. If $E$ is any equivalence relation on
$\kappa^{+}$, then there exists a notion of forcing
$\mathbb{P}$ such that
\begin{enumerate}
\item[(a)] $\mathbb{P}$ preserves cofinalities and cardinalities;
\item[(b)] $\mathbb{P}$ does not adjoin any new $\kappa$-sequences 
of ordinals;
\item[(c)] each graph $\Gamma_{\alpha}$ remains rigid in $V^{\mathbb{P}}$;
\item[(d)] $\Gamma_{\alpha} \simeq \Gamma_{\beta}$ in $V^{\mathbb{P}}$ if{f}
$\alpha \mathrel{E} \beta$. 
\end{enumerate}
\end{Hyp}

Our argument will use the normaliser
tower technique, which was introduced in \cite{t1}.

\begin{Def} \label{D:norm}
If $G$ is a subgroup of the group $H$, then the normaliser tower
of $G$ in $H$ is defined inductively as follows.
\begin{enumerate}
\item[(a)] $N_{0}(G) = G$.
\item[(b)] If $\alpha = \beta + 1$, then
$N_{\alpha}(G) = N_{H} \left( N_{\beta}(G) \right)$.
\item[(c)] If $\alpha$ is a limit ordinal, then
$N_{\alpha}(G) = \underset{\beta < \alpha}{\bigcup}N_{\beta}(G)$.
\end{enumerate}
\end{Def}

The definition of the normaliser tower is motivated by the following
observation.

\begin{Prop} \label{P:norm}
Let $G$ be a centreless group, and let
$\left( G_{\alpha} \mid \alpha \leq \tau(G) \right)$ be the automorphism
tower of $G$. Then for each $\alpha \leq \tau(G)$,
$G_{\alpha} = N_{\alpha}(G)$, where $N_{\alpha}(G)$ is the 
$\alpha^{th}$ group in the normaliser tower of $G$ in
$G_{\tau(G)}$.
\end{Prop}

\begin{proof}
Let $\gamma = \tau(G)$. We will show that
$N_{G_{\gamma}}(G_{\alpha}) = G_{\alpha + 1}$
for all $\alpha < \gamma$. Since the inclusion
$G_{\alpha} \leqslant G_{\alpha +1}$ is isomorphic to
the inclusion $\Inn (G_{\alpha}) \leqslant \Aut (G_{\alpha})$,
it follows that $G_{\alpha + 1} \leqslant N_{G_{\gamma}}(G_{\alpha})$.
Conversely, suppose that $g \in N_{G_{\gamma}}(G_{\alpha})$.
Then $g$ induces an automorphism of $G_{\alpha}$ via conjugation.
Hence there exists $h \in G_{\alpha + 1}$ such that
$h x h^{-1} = g x g^{-1}$ for all $x \in G_{\alpha}$. Thus
$h^{-1}g \in C_{G_{\gamma}}(G_{\alpha})$. By Lemma 8.1.1
of Hulse \cite{h}, $ C_{G_{\gamma}}(G_{\alpha}) = 1$. Hence
$g = h \in G_{\alpha +1}$. Consequently,
$N_{G_{\gamma}}(G_{\alpha}) = G_{\alpha + 1}$.
\end{proof}

The following lemma, which was essentially proved in \cite{t1},
will enable us to convert  normaliser towers into 
corresponding automorphism towers. (The proof makes use of the
assumption that $PSL (2,K)$ is simple. This is true if and only if
$|K| > 3$.)

\begin{Lem} \label{L:psl}
Let $K$ be a field such that $|K| > 3$ and let $H$ be a subgroup of $\Aut(K)$.
Let
\[
G = PGL(2,K) \rtimes H \leqslant P \varGamma L(2,K) =
PGL(2,K) \rtimes \Aut(K).
\]
Then $G$ is a centreless group; and for each $\alpha$,
$G_{\alpha} = PGL(2,K) \rtimes N_{\alpha}(H)$, where
$N_{\alpha}(H)$ is the $\alpha^{th}$ group in the normaliser
tower of $H$ in $\Aut(K)$. 
\end{Lem}
\begin{flushright}
$\square$
\end{flushright}

We will also make use of the following result.

\begin{Lem}[Fried and Koll\'{a}r \cite{f}] \label{L:field}
Let $\Gamma = \langle X, E \rangle$ be any graph. Then there exists a
field $K_{\Gamma}$ which satisfies the following conditions.
\begin{enumerate}
\item[(a)] $X \subseteq K_{\Gamma}$.
\item[(b)] If $\mathbb{P}$ is a (possibly trivial) notion of forcing and
$M = V^{\mathbb{P}}$, then
\begin{enumerate}
\item[(i)] $\pi [X] = X$ for all $\pi \in \Aut ^{M} (K_{\Gamma})$; and
\item[(ii)] the restriction mapping, $\pi \mapsto \pi \res X$ is an 
isomorphism from $\Aut ^{M} (K_{\Gamma})$ onto $\Aut ^{M} (\Gamma)$.
\end{enumerate}
\end{enumerate}
\end{Lem}

\begin{proof}
This follows from the observation that the construction of Fried and
Koll\'{a}r in \cite{f} is upwards absolute.
\end{proof}

Most of our effort in this section will go into proving the following
result.

\begin{Thm} \label{T:norm}
Suppose that Hypothesis \ref{H:graph} holds. Then for every infinite
cardinal $\lambda$ and every ordinal $\alpha < \lambda$, there exist a
graph $\Gamma$ and a subgroup $H \leqslant \Aut( \Gamma)$ with
the following properties.
\begin{enumerate}
\item[(a)] The normaliser tower of $H$ in $\Aut(\Gamma)$ terminates
in exactly $\alpha$ steps.
\item[(b)] If $\beta$ is any ordinal such that
$1 \leq \beta < \lambda$, then there exists a notion of
forcing $\mathbb{P}$, which preserves cofinalities and cardinalities,
such that the normaliser tower of $H$ in 
$\Aut^{V^{\mathbb{P}}}(\Gamma)$ terminates in exactly $\beta$
steps.
\end{enumerate}
\end{Thm}

\begin{Cor} \label{C:norm}
Suppose that Hypothesis \ref{H:graph} holds. Then for every infinite
cardinal $\lambda$ and every ordinal $\alpha < \lambda$, there exists a
centreless group $G$ with the following properties. 
\begin{enumerate}
\item[(a)] $\tau(G) = \alpha$. 
\item[(b)] If $\beta$ is any ordinal such that
$1 \leq \beta < \lambda$, then there exists a notion of
forcing $\mathbb{P}$, which preserves cofinalities and cardinalities,
such that $\tau^{V^{\mathbb{P}}}(G) = \beta$.
\end{enumerate}
\end{Cor}

\begin{proof}
Let $\Gamma$ and $H \leq \Aut(\Gamma)$ be the graph and subgroup
which are given by Theorem \ref{T:norm}. Let $K_{\Gamma}$ be the
corresponding field which is given by Lemma \ref{L:field}. By Lemma
\ref{L:psl}, the centreless group
$ G = PGL(2,K_{\Gamma}) \rtimes H$ satisfies our requirements.
\end{proof}

Now we will begin the proof of Theorem \ref{T:norm}.

\begin{Def} \label{D:product}
\begin{enumerate}
\item[(a)] Suppose that $\Gamma_{i} = \left( X_{i}, E_{i} \right)$
is a graph for each $i \in I$. Then the {\em direct sum\/} of the
graphs $\{ \Gamma_{i} \mid i \in I \}$ is defined to be the graph
\[
\bigoplus_{i \in I} \Gamma_{i} = \left( \bigsqcup_{i \in I} X_{i},
\bigsqcup_{i \in I} E_{i} \right).
\]
Here $\bigsqcup_{i \in I} X_{i}$ and $\bigsqcup_{i \in I} E_{i}$
denote the {\em disjoint unions\/} of the sets of vertices and edges
respectively.
\item[(b)] Suppose further that $H_{i} \leqslant \Aut(\Gamma_{i})$
for each $i \in I$. Then the {\em direct product\/} of the permutation
groups $\{ \left( H_{i}, \Gamma_{i} \right) \mid i \in I \}$ is defined
to be the permutation group
\[
\prod_{i \in I} \left( H_{i}, \Gamma_{i} \right) =
\left( \prod_{i \in I} H_{i}, \bigoplus_{i \in I} \Gamma_{i} \right),
\]
where $\prod_{i \in I} H_{i}$ acts on
$\bigoplus_{i \in I} \Gamma_{i}$ in the obvious manner. If
$I = \{ 1, 2 \}$, then we write
$\prod_{i \in I} \left( H_{i}, \Gamma_{i} \right) =
\left( H_{1}, \Gamma_{1} \right) \times 
\left( H_{2}, \Gamma_{2} \right) =
\left( H_{1} \times H_{2}, \Gamma_{1} \oplus \Gamma_{2} \right)$.
\end{enumerate}
\end{Def}

\begin{Def} \label{D:ind}
Let $\Gamma$ be a rigid connected graph. For each $\alpha$, we define 
permutation groups $\left( H_{\alpha}( \Gamma), \mathcal{G}_{\alpha}(\Gamma) \right)$
and $\left( F_{\alpha}( \Gamma), \mathcal{G}_{\alpha}(\Gamma)  \right)$
inductively as follows.
\begin{enumerate}
\item[(a)] $\left( H_{0}( \Gamma), \mathcal{G}_{0}(\Gamma)  \right) =
\left( F_{0}( \Gamma), \mathcal{G}_{0}(\Gamma)  \right) = 
\left( \Aut (\Gamma), \Gamma \right) =
\left( 1, \Gamma \right)$.
\item[(b)] If $\alpha > 0$, then we define
\[
\left( H_{\alpha}(\Gamma) , \mathcal{G}_{\alpha}(\Gamma)  \right) =
\left( F_{0}(\Gamma), \mathcal{G}_{0}(\Gamma)  \right) \times
\prod_{\beta < \alpha} \left( F_{\beta}(\Gamma),
\mathcal{G}_{\beta}(\Gamma)  \right) ,
\]
and we define
$F_{\alpha}(\Gamma)$ to be the terminal group of the normaliser
tower of $H_{\alpha}(\Gamma)$ in $\Aut(\mathcal{G}_{\alpha}(\Gamma))$.
\end{enumerate}
\end{Def}

In the proof of the following lemma, we will need to study the blocks of
imprimitivity of $F_{\alpha}(\Gamma)$ in its action on the set of
connected components of $\mathcal{G}_{\alpha}(\Gamma)$. Recall that
if $\left( G, \Omega \right)$ is a transitive permutation group, then the 
nonempty subset $Z$ of $\Omega$ is a {\em block of imprimitivity\/} if
for each $g \in G$, either $g[Z] = Z$ or
$g[Z] \cap Z = \emptyset$. In this case, we obtain a $G$-invariant
equivalence relation $E$ on $\Omega$ corresponding to the partition
$\{ g[Z] \mid g \in G \}$.

\begin{Lem} \label{L:stops}
If $\Gamma$ is a rigid connected graph,
then normaliser tower of $H_{\alpha}(\Gamma)$ in
$\Aut(\mathcal{G}_{\alpha}(\Gamma))$ terminates in exactly $\alpha$ steps.
\end{Lem}

\begin{proof}
During the course of this proof, we will need a more explicit 
definition of the graph $\mathcal{G}_{\alpha}(\Gamma)$. So for
each ordinal $\alpha$, we define the graphs 
$\mathcal{G}_{\alpha}(\Gamma)$ and $\mathcal{G}_{\alpha}^{1}(\Gamma)$ 
inductively as follows.
\begin{enumerate}
\item[(a)] $\mathcal{G}_{0}(\Gamma) = \Gamma$.
\item[(b)] Suppose that $\alpha = \beta +1$ and that
$\mathcal{G}_{\beta}(\Gamma)$ has been defined. Let
$\mathcal{G}_{\beta}^{1}(\Gamma)$ be a graph such that
$\mathcal{G}_{\beta}^{1}(\Gamma) \cap
\mathcal{G}_{\beta}(\Gamma) = \emptyset$ and
$\mathcal{G}_{\beta}^{1}(\Gamma) \simeq
\mathcal{G}_{\beta}(\Gamma)$. Then we define
$\mathcal{G}_{\beta +1}(\Gamma) =
\mathcal{G}_{\beta}(\Gamma) \cup
\mathcal{G}_{\beta}^{1}(\Gamma)$.
\item[(c)] If $\alpha$ is a limit ordinal, then we
define $\mathcal{G}_{\alpha}(\Gamma) =
\underset{\beta < \alpha}{\bigcup}
\mathcal{G}_{\beta}(\Gamma)$.
\end{enumerate}
In particular, for each ordinal $\beta$, we have that
\[
\left( H_{\beta +1}(\Gamma), \mathcal{G}_{\beta +1}(\Gamma) \right)
=
\left( H_{\beta}(\Gamma), \mathcal{G}_{\beta}(\Gamma) \right)
\times
\left( F_{\beta}(\Gamma), \mathcal{G}^{1}_{\beta}(\Gamma) \right).
\]
For each ordinal $\alpha$, let $\Delta_{\alpha}$ and
$\Delta_{\alpha}^{1}$ be the sets of connected components of
the graphs $\mathcal{G}_{\alpha}(\Gamma)$ and
$\mathcal{G}_{\alpha}^{1}(\Gamma)$ respectively. Then for all
$\beta < \alpha$, we have that
\[
\Delta_{\alpha} = \Delta_{\beta} \cup 
\bigcup_{\beta \leq \gamma < \alpha} \Delta^{1}_{\gamma}.
\]
Since $\Gamma$ is a rigid 
connected graph, $\Aut ( \mathcal{G}_{\alpha}(\Gamma) )$ can be
identified naturally with $\Sym ( \Delta_{\alpha} )$. This allows us to
regard $H_{\alpha}( \Gamma)$ and $F_{\alpha}( \Gamma)$ as
subgroups of $\Sym ( \Delta_{\alpha} )$.
We will need to consider direct products of permutation groups of the form
$\left( F_{\alpha}( \Gamma ), \Delta_{\alpha} \right)$ or
$\left( F_{\alpha}( \Gamma ), \Delta^{1}_{\alpha} \right)$. We will regard each
set $\Delta_{\alpha}$ and $\Delta^{1}_{\alpha}$
as a null graph, and continue to use the notation
introduced in Definition \ref{D:product}. Thus
\[
\left( H_{\alpha}(\Gamma) , \Delta_{\alpha} \right) =
\left( F_{0}(\Gamma), \Delta_{0} \right) \times
\prod_{\beta < \alpha} \left( F_{\beta}(\Gamma), \Delta^{1}_{\beta} \right) .
\]

We will prove the following statements by a simultaneous induction on
$\alpha \geq 0$.
\begin{enumerate}
\item[$(1_{\alpha})$] $F_{\alpha}(\Gamma)$ acts transitively on
$\Delta_{\alpha}$.
\item[$(2_{\alpha})$] Let $\Delta_{0} = \{ v_{0} \}$. Then
$\{ \Delta_{\beta} \mid \beta \leq \alpha \}$ is the set of blocks $Z$
of imprimitivity of $\left( F_{\alpha}( \Gamma ), \Delta_{\alpha} \right)$
such that $v_{0} \in Z$.
\item[$(3_{\alpha})$] For each $\beta \leq \alpha$, let $E_{\beta}$ be the
$F_{\alpha}(\Gamma)$-invariant equivalence relation corresponding to the 
partition $\{ g \left[ \Delta_{\beta} \right] \mid g \in F_{\alpha}(\Gamma) \}$.
Then for each $\beta \leq \gamma < \alpha$, the set $\Delta^{1}_{\gamma}$
is a union of $E_{\beta}$-equivalence classes. (Strictly speaking, we should
write $E^{\alpha}_{\beta}$ to indicate that this is the 
$F_{\alpha}(\Gamma)$-invariant equivalence relation
on $\Delta_{\alpha}$ corresponding to the block of imprimitivity $\Delta_{\beta}$.
However, this slight abuse of notation should not cause any confusion.)
\item[$(4_{\alpha})$] If $\beta < \alpha$, then
\[
\left( N_{\beta}( H_{\alpha} ( \Gamma )), \Delta_{\alpha} \right)
= \left( F_{\beta}( \Gamma ), \Delta_{\beta} \right) \times
\prod_{\beta \leq \gamma < \alpha}
\left( F_{\gamma}( \Gamma ), \Delta^{1}_{\gamma} \right).
\]
\item[$(5_{\alpha})$] Furthermore, if $\beta < \alpha$, then
$ N_{\beta}( H_{\alpha} ( \Gamma ))$ is the stabiliser of the partition
$\{ \Delta_{\beta} \} \cup \{ \Delta^{1}_{\gamma} \mid
\beta \leq \gamma < \alpha \}$ in $F_{\alpha}(\Gamma)$.
\item[$(6_{\alpha})$] $N_{\alpha}( H_{\alpha} ( \Gamma ))$ is
self-normalising in $\Sym ( \Delta_{\alpha} )$; and so
$N_{\alpha}( H_{\alpha} ( \Gamma )) = F_{\alpha}( \Gamma )$.
\end{enumerate}

Let $\beta \leq \alpha$ and let $v$ be any point of $\Delta_{\beta}$. Let
$\mathcal{B}_{\beta}(v)$ be the set of blocks $Z$ of imprimitivity of
$\left( F_{\beta}( \Gamma ), \Delta_{\beta} \right)$ such that $v \in Z$.
Then conditions $(1_{\beta})$ and $(2_{\beta})$ imply that
$\mathcal{B}_{\beta}(v)$, ordered by inclusion, is a well-ordering
of order-type $\beta +1$. Hence conditions 
$(1_{\beta})$ and $(2_{\beta})$ for all $\beta \leq \alpha$ yield the
following statement.
\begin{enumerate}
\item[$(7_{\alpha})$] If $\beta < \gamma \leq \alpha$, then
$\left( F_{\beta}( \Gamma ), \Delta_{\beta} \right)$ and
$\left( F_{\gamma}( \Gamma ), \Delta_{\gamma} \right)$ are
nonisomorphic permutation groups.
\end{enumerate}

Note that $\left| \Delta_{n} \right| = 2^{n}$ for all $n \in \omega$.
It is easily checked that the result holds for $\alpha = 0,1,2$.
Furthermore, $F_{0}(\Gamma) = 1$, 
$F_{1}(\Gamma) = \Sym(2)$ and
$F_{2}(\Gamma) = \Sym(2) \Wr \Sym(2)$. In the successor stage
of the argument, we will appeal to a result of Neumann \cite{n}
on the automorphism groups of wreath products
$A \Wr B$. The hypotheses of this result require that $A$ should
not be a ``special dihedral'' group. For our purposes, it is enough to
know that if $\Sym(2) \Wr \Sym(2) \leqslant A$, then $A$ is not a
``special dihedral'' group.

Now suppose that $\alpha \geq 2$ and that the result holds for all
$\beta \leq \alpha$. Remember that if 
$\left( G, \Omega \right)$ is a permutation group and 
$\pi \in \Sym( \Omega )$ normalises $G$, then $\pi$ permutes the
orbits of $G$. Furthermore, if $X$ and $Y$ are $G$-orbits and
$\pi[X] = Y$, then $G$ must induce isomorphic permutation groups
via its actions on $X$ and $Y$.
Using this observation, together with the inductive
hypotheses for $\beta \leq \alpha$, we see that
\[
\left( N_{\beta}( H_{\alpha +1} ( \Gamma )), \Delta_{\alpha +1} \right)
= \left( F_{\beta}( \Gamma ), \Delta_{\beta} \right) \times
\prod_{\beta \leq \gamma < \alpha +1}
\left( F_{\gamma}( \Gamma ), \Delta^{1}_{\gamma} \right)
\]
for all $\beta \leq \alpha$. In particular, we have that
\[
\left( N_{\alpha}( H_{\alpha +1} ( \Gamma )), \Delta_{\alpha +1} \right)
= \left( F_{\alpha}( \Gamma ), \Delta_{\alpha} \right) \times
\left( F_{\alpha}( \Gamma ), \Delta^{1}_{\alpha} \right).
\]
It follows easily that
\[
\left( N_{\alpha +1}( H_{\alpha +1} ( \Gamma )), \Delta_{\alpha +1} \right)
= \left( F_{\alpha}(\Gamma) \Wr \Sym(2) ,
\Delta_{\alpha} \sqcup \Delta^{1}_{\alpha} \right).
\]
Notice that we have already established conditions $(1_{\alpha +1})$ and
$(4_{\alpha +1})$. Next suppose that $g \in \Sym (\Delta_{\alpha +1})$
normalises $N_{\alpha +1}( H_{\alpha +1} ( \Gamma )) =
F_{\alpha}(\Gamma) \Wr \Sym(2)$. By Theorem 9.12 of Neumann \cite{n},
the base group $F_{\alpha}(\Gamma) \times F_{\alpha}(\Gamma)$ 
is a characteristic subgroup of the
wreath product $ F_{\alpha}(\Gamma) \Wr \Sym(2)$. Hence $g$ must also
normalise $N_{\alpha}(H_{\alpha +1}(\Gamma)) =
F_{\alpha}(\Gamma) \times F_{\alpha}(\Gamma)$; and so
$g \in N_{\alpha +1}( H_{\alpha +1} ( \Gamma ))$. Thus condition
$(6_{\alpha +1})$ also holds. It is now easily checked that conditions
$(2_{\alpha +1})$, $(3_{\alpha +1})$ and $(5_{\alpha +1})$ hold. Thus the result
holds for $\alpha +1$.

Now suppose that $\lambda$ is a limit ordinal, and that the result holds for
all $\alpha < \lambda$. Once again, it is easy to see that conditions
$(4_{\lambda})$ and $(1_{\lambda})$ hold. It is also easily checked that
the following statements hold.
\begin{enumerate}
\item[$(2_{\lambda})^{\prime}$] For each $\beta \leq \lambda$, the set
$\Delta_{\beta}$ is a block of imprimitivity of
$\left( N_{\lambda}(H_{\lambda}(\Gamma)), \Delta_{\lambda} \right)$.
\item[$(3_{\lambda})^{\prime}$] For each $\beta \leq \lambda$, let
$E_{\beta}$ be the $N_{\lambda}(H_{\lambda}(\Gamma))$-invariant
equivalence relation corresponding to the partition
$\{ g \left[ \Delta_{\beta} \right] \mid g \in N_{\lambda}(H_{\lambda}(\Gamma)) \}$.
Then for each $\beta \leq < \lambda$, the set $\Delta^{1}_{\gamma}$ is a union
of $E_{\beta}$-equivalence classes.
\item[$(5_{\lambda})^{\prime}$] If $\beta < \lambda$, then
$N_{\beta}( H_{\lambda} ( \Gamma ))$ is the stabiliser of the partition
$\{ \Delta_{\beta} \} \cup \{ \Delta^{1}_{\gamma} \mid
\beta \leq \gamma < \lambda \}$ in $N_{\lambda}(H_{\lambda}(\Gamma))$.
\end{enumerate}

Thus it is enough to prove the following two claims.

\begin{Claim} \label{C:block}
If $Z$ is a block of imprimitivity of
$\left( N_{\lambda}(H_{\lambda}(\Gamma)), \Delta_{\lambda} \right)$
such that $v_{0} \in Z$, then $Z = \Delta_{\beta}$ for some 
$\beta \leq \lambda$.
\end{Claim}

\begin{Claim} \label{C:end}
$N_{\lambda}( H_{\lambda} ( \Gamma ))$ is
self-normalising in $\Sym ( \Delta_{\lambda} )$; and so
$N_{\lambda}( H_{\lambda} ( \Gamma )) = F_{\lambda}( \Gamma )$.
\end{Claim}

\begin{proof}[Proof of Claim \ref{C:block}]
If there exists $\gamma < \lambda$ such that $Z \subseteq \Delta_{\beta}$,
then it follows from the inductive hypotheses that
$Z = \Delta_{\beta}$ for some $\beta \leq \gamma$. Hence we can suppose
that the set
$I = \{ \gamma < \lambda \mid Z \cap \Delta^{1}_{\gamma} \neq 
\emptyset \}$ is cofinal in $\lambda$. Fix some $\gamma \in I$.
By condition $(4_{\lambda})$,
\[
\left( N_{\gamma}( H_{\lambda} ( \Gamma )), \Delta_{\lambda} \right)
= \left( F_{\gamma}( \Gamma ), \Delta_{\gamma} \right) \times
\prod_{\gamma \leq \xi < \lambda}
\left( F_{\xi}( \Gamma ), \Delta^{1}_{\xi} \right).
\]
Hence for each $x \in \Delta_{\gamma}$, there exists an element
$g \in N_{\gamma}( H_{\lambda} ( \Gamma )) \leqslant
N_{\lambda}( H_{\lambda} ( \Gamma ))$ such that
\begin{enumerate}
\item[(i)] $g(v_{0}) = x$; and
\item[(ii)] $g(y) =y$ for all $y \in \Delta^{1}_{\gamma}$.
\end{enumerate}
By condition (ii), $g[Z] \cap Z \neq \emptyset$ and hence
$x \in g[Z] = Z$. Thus $\Delta_{\gamma} \subseteq Z$ for
each $\gamma \in I$, and so $Z = \Delta_{\lambda}$.
\end{proof}

\begin{proof}[Proof of Claim \ref{C:end}]
Note that for each $\pi \in N_{\lambda}( H_{\lambda} ( \Gamma ))$,
there exists $\beta < \lambda$ such that
$\pi \in N_{\beta}( H_{\lambda} ( \Gamma ))$; and so
$\pi \left[ \Delta^{1}_{\gamma} \right] = \Delta^{1}_{\gamma}$
for all $\beta \leq \gamma < \lambda$.
Suppose that $g \in \Sym( \Delta_{\lambda})$ normalises
$N_{\lambda}(H_{\lambda}(\Gamma))$. Then $g$ must permute
the set $\{ E_{\beta} \mid \beta \leq \lambda \}$ of
$N_{\lambda}(H_{\lambda}(\Gamma))$-invariant equivalence relations
on $\Delta_{\lambda}$. Since the set 
$\{ E_{\beta} \mid \beta \leq \lambda \}$ is well-ordered under
inclusion, it follows that $E_{\beta}$ is also $g$-invariant for
each $\beta \leq \lambda$.

Next we will show that
there exists $\beta < \lambda$ such that $g \left[ \Delta^{1}_{\gamma} \right]
= \Delta^{1}_{\gamma}$ for all $\beta \leq \gamma < \lambda$.
For each $\gamma < \lambda$, let 
$C^{1}_{\gamma} = g \left[ \Delta^{1}_{\gamma} \right]$. Suppose that there
exists a cofinal subset $I \subseteq \lambda$ such that 
$C^{1}_{\gamma} \neq \Delta^{1}_{\gamma}$ for all $\gamma \in I$.
Fix some $\gamma \in I$. Since $\Delta^{1}_{\gamma}$ is an
$E_{\gamma}$-equivalence class, it follows that
$C^{1}_{\gamma} = g \left[ \Delta^{1}_{\gamma} \right]$ is also an
$E_{\gamma}$-equivalence class. Using condition $(3_{\lambda})^{\prime}$,
it follows that either
\begin{enumerate}
\item[(i)] $C^{1}_{\gamma} = \Delta_{\gamma}$, or
\item[(ii)] there exists $f(\gamma) > \gamma$ such that
$C^{1}_{\gamma} \varsubsetneq \Delta^{1}_{f(\gamma)}$.
\end{enumerate}
Clearly we can assume that condition (ii) holds for all 
$\gamma \in I$. Furthermore, by passing to a suitable subset of $I$ if
necessary, we can assume that the resulting function
$f : I \to \lambda$ is injective. Remember that
\[
\left( H_{\lambda} ( \Gamma ), \Delta_{\lambda} \right)
= \left( F_{0}( \Gamma ), \Delta_{0} \right) \times
\prod_{0 \leq \xi < \lambda}
\left( F_{\xi}( \Gamma ), \Delta^{1}_{\xi} \right).
\]
Since each $F_{\xi}( \Gamma )$ acts transitively on
$\Delta^{1}_{\xi}$, there exists an element
$\psi \in H_{\lambda}(\Gamma) \leqslant N_{\lambda}(H_{\lambda}(\Gamma))$ 
such that $\psi \left[ C^{1}_{\gamma} \right] \neq C^{1}_{\gamma}$
for all $\gamma \in I$. Let
$\pi = g^{-1} \psi g \in N_{\lambda}(H_{\lambda}(\Gamma))$. Then
$\pi \left[ \Delta^{1}_{\gamma} \right] \neq \Delta^{1}_{\gamma}$
for all $\gamma \in I$, which is a contradiction.

Thus there exists $\beta < \lambda$ such that
$g \in \Sym( \Delta_{\beta}) \times 
\prod_{\beta \leq \gamma < \lambda} \Sym( \Delta^{1}_{\gamma})$.
Since $N_{\beta}( H_{\lambda} ( \Gamma ))$ is the stabiliser of the partition
$\{ \Delta_{\beta} \} \cup \{ \Delta^{1}_{\gamma} \mid
\beta \leq \gamma < \lambda \}$ in $N_{\lambda}(H_{\lambda}(\Gamma))$,
it follows that $g$ normalises $ N_{\beta}( H_{\lambda} ( \Gamma ))$.
Thus $g \in N_{\beta +1}( H_{\lambda} ( \Gamma )) \leqslant
N_{\lambda}(H_{\lambda}(\Gamma))$.
\end{proof}

This completes the proof of Lemma \ref{L:stops}.
\end{proof}

\begin{Def} \label{D:down}
Let $\Gamma$ be a connected rigid graph. If
$1 \leq \beta < \alpha$, then we define
\[
\left( D^{\alpha}_{\beta}(\Gamma), \mathcal{G}^{\alpha}_{\beta}(\Gamma)  \right)
= \left( H_{\alpha}( \Gamma ), \mathcal{G}_{\alpha}(\Gamma)  \right) \times
\left( F_{\beta} (\Gamma), \mathcal{G}_{\beta}(\Gamma)  \right) \times
\left( F_{\beta}(\Gamma), \mathcal{G}_{\beta}(\Gamma)  \right).
\]
\end{Def}

\begin{Lem} \label{L:down}
If $\Gamma$ is a connected rigid graph, then the normaliser tower
of $ D^{\alpha}_{\beta}(\Gamma)$ in $\Aut(\mathcal{G}^{\alpha}_{\beta}(\Gamma))$
terminates in exactly $\beta$ steps.
\end{Lem}

\begin{proof}
Let
$\left( B, \Gamma^{\prime} \right) =
\left( F_{\beta} (\Gamma), \mathcal{G}_{\beta}(\Gamma)  \right) \times
\left( F_{\beta} (\Gamma), \mathcal{G}_{\beta}(\Gamma)  \right) \times
\left( F_{\beta} (\Gamma), \mathcal{G}_{\beta}(\Gamma)  \right)$, 
and let $ F_{\beta} (\Gamma) \mathbin{\Wr} \Sym (3) =
B \rtimes \Sym (3)$ be the associated wreath product. By rearranging
the order of its factors, we can identify 
$\left( D^{\alpha}_{\beta}(\Gamma), \mathcal{G}^{\alpha}_{\beta}(\Gamma)  \right)$
with
\[
\left( H_{\beta}(\Gamma), \mathcal{G}_{\beta}(\Gamma) \right) \times
\left( B, \Gamma^{\prime} \right) \times
\prod_{\beta < \gamma < \alpha} 
\left( F_{\gamma} (\Gamma), \mathcal{G}_{\gamma}(\Gamma)  \right).
\]
Arguing as in the proof of Lemma \ref{L:stops}, we find that the
$\beta^{th}$ element of the normaliser tower of
$ D^{\alpha}_{\beta}(\Gamma)$ in $\Aut(\mathcal{G}^{\alpha}_{\beta}(\Gamma))$
is
\[
\left( F_{\beta}(\Gamma), \mathcal{G}_{\beta}(\Gamma) \right) \times
\left(F_{\beta} (\Gamma) \mathbin{\Wr} \Sym (3), \Gamma^{\prime} \right) \times
\prod_{\beta < \gamma < \alpha} 
\left( F_{\gamma} (\Gamma), \mathcal{G}_{\gamma}(\Gamma)  \right);
\]
and also that this group is self-normalising in 
$\Aut(\mathcal{G}^{\alpha}_{\beta}(\Gamma))$.
\end{proof}

Now let $\lambda$ be any infinite cardinal and let $\alpha$ be any ordinal
such that $\alpha < \lambda$. Choose a regular cardinal $\kappa$ such that
$\lambda \leq \kappa$. Let 
$\{ \Gamma_{\gamma} \mid \gamma < \kappa^{+} \}$ be the set of
pairwise nonisomorphic connected rigid graphs given by Hypothesis
\ref{H:graph}. If $\alpha \geq 1$, then we define
\[
\left( B_{\alpha}, \Gamma^{\alpha} \right) = 
\prod_{1 \leq \beta < \alpha} \Bigl( 
\left( F_{\beta }(\Gamma_{\beta}), \mathcal{G}_{\beta }(\Gamma_{\beta}) \right)
\times
\left( F_{\beta }(\Gamma_{\beta}), \mathcal{G}_{\beta }(\Gamma_{\beta}) \right)
\Bigr)
\]
and
\[
\left( H, \Gamma \right) =
\left( B_{\alpha}, \Gamma^{\alpha} \right) \times
\left( H_{\alpha}(\Gamma_{\alpha}), \mathcal{G}_{\alpha}(\Gamma_{\alpha}) \right)
\times
\prod_{\alpha \leq \gamma < \lambda} 
\left( F_{\gamma}(\Gamma_{\gamma +1}), \mathcal{G}_{\gamma}(\Gamma_{\gamma +1}) \right).
\]
If $\alpha =0$, then we define
\[
\left( H, \Gamma \right) =
\left( F_{0}(\Gamma_{0}), \mathcal{G}_{0}(\Gamma_{0}) \right)
\times
\prod_{0 \leq \gamma < \lambda} 
\left( F_{\gamma}(\Gamma_{\gamma +1}), \mathcal{G}_{\gamma}(\Gamma_{\gamma +1}) \right).
\]
We will show that $\Gamma$ and $H \leqslant \Aut(\Gamma)$ satisfy the requirements of
Theorem \ref{T:norm}.

\begin{Lem} \label{L:almost}
The normaliser tower of $H$ in $\Aut(\Gamma)$ terminates in exactly
$\alpha$ steps.
\end{Lem}

\begin{proof}
For example, suppose that $\alpha \geq 1$. Then
\[
\Aut ( \Gamma) = \Aut ( \Gamma^{\alpha}) \times
\Aut ( \mathcal{G}_{\alpha}( \Gamma_{\alpha})) \times
\prod_{\alpha \leq \gamma < \lambda}
\Aut ( \mathcal{G}_{\gamma}( \Gamma_{\gamma +1})).
\]
Using Lemma \ref{L:stops}, we see that
\begin{enumerate}
\item[(a)] the normaliser tower of $B_{\alpha}$ in $\Aut ( \Gamma^{\alpha})$
terminates in exactly 1 step;
\item[(b)] the normaliser tower of $H_{\alpha}( \Gamma_{\alpha})$ in
$\Aut ( \mathcal{G}_{\alpha}( \Gamma_{\alpha}))$ terminates in exactly
$\alpha$ steps; and
\item[(c)] the normaliser tower of
$\prod_{\alpha \leq \gamma < \lambda} F_{\gamma}( \Gamma_{\gamma +1})$
in $\prod_{\alpha \leq \gamma < \lambda}
\Aut ( \mathcal{G}_{\gamma}( \Gamma_{\gamma +1}))$ terminates in exactly
0 steps.
\end{enumerate}
It follows that the normaliser tower of $H$ in $\Aut(\Gamma)$ terminates in exactly
$\alpha$ steps.
\end{proof}

\begin{Lem} \label{L:there}
If $\beta$ is any ordinal such that
$1 \leq \beta < \lambda$, then there exists a notion of
forcing $\mathbb{P}$, which preserves cofinalities and cardinalities,
such that the normaliser tower of $H$ in 
$\Aut^{V^{\mathbb{P}}}(\Gamma)$ terminates in exactly $\beta$
steps.
\end{Lem}

\begin{proof}
Clearly we can suppose that $\beta \neq \alpha$. There are two cases to consider.
First suppose that $1 \leq \beta < \alpha$. Let $E$ be the equivalence relation
on $\kappa^{+}$ such that
\[
\gamma \mathrel{E} \delta \text{ if{f} } \{ \gamma, \delta \} = \{ \alpha, \beta \}
\text{ or } \gamma = \delta ;
\]
and let $\mathbb{P}$ be the corresponding notion of forcing, given by
Hypothesis \ref{H:graph}. Using the facts that
\begin{enumerate}
\item[(a)] each graph $\Gamma_{\delta}$ remains rigid in $V^{\mathbb{P}}$, and
\item[(b)] $\mathbb{P}$ does not adjoin any new $\kappa$-sequences of
ordinals,
\end{enumerate}
we see that
$\left( H_{\gamma}(\Gamma_{\delta}), 
\mathcal{G}_{\gamma}(\Gamma_{\delta}) \right) ^{V^{\mathbb{P}}}
= \left( H_{\gamma}(\Gamma_{\delta}), 
\mathcal{G}_{\gamma}(\Gamma_{\delta}) \right)$ and
$\left( F_{\gamma}(\Gamma_{\delta}), 
\mathcal{G}_{\gamma}(\Gamma_{\delta}) \right) ^{V^{\mathbb{P}}}
= \left( F_{\gamma}(\Gamma_{\delta}), 
\mathcal{G}_{\gamma}(\Gamma_{\delta}) \right)$ for all
$\gamma$, $\delta < \kappa^{+}$. Let
\[
\left( B^{\prime}, \Gamma^{\prime} \right) =
\prod_{ \substack{ 1 \leq \gamma < \alpha \\
\gamma \neq \beta }} \left( 
\left( F_{\gamma }(\Gamma_{\gamma}), \mathcal{G}_{\gamma }(\Gamma_{\gamma}) \right)
\times
\left( F_{\gamma }(\Gamma_{\gamma}), \mathcal{G}_{\gamma }(\Gamma_{\gamma}) \right)
\right).
\]
Then in $V^{\mathbb{P}}$, $\left( H, \Gamma \right)$ is isomorphic to
\[
\left( B^{\prime}, \Gamma^{\prime} \right) \times
\left( D^{\alpha}_{\beta}(\Gamma_{\alpha}), \mathcal{G}^{\alpha}_{\beta}(\Gamma_{\alpha}) \right)
\times
\prod_{\alpha \leq \gamma < \lambda} 
\left( F_{\gamma}(\Gamma_{\gamma +1}), \mathcal{G}_{\gamma}(\Gamma_{\gamma +1}) \right).
\]
Hence the normaliser tower of $H$ in 
$\Aut^{V^{\mathbb{P}}}(\Gamma)$ terminates in exactly $\beta$
steps.

Now suppose that $\alpha < \beta < \lambda$. We will only deal with the
case when $\alpha \geq 1$. (The case when $\alpha = 0$ is almost
identical.) Now let $E$ be the equivalence relation on $\kappa^{+}$
such that 
\[
\gamma \mathrel{E} \delta \text{ if{f} }
\alpha \leq \gamma , \delta < \beta +1 \text{ or } \gamma = \delta ;
\]
and let $\mathbb{P}$ be the corresponding notion of forcing, given by
Hypothesis \ref{H:graph}. Then 
in $V^{\mathbb{P}}$, $\left( H, \Gamma \right)$ is isomorphic to
\[
\left( B_{\alpha}, \Gamma^{\alpha} \right) \times
\left( H_{\beta}(\Gamma_{\alpha}), \mathcal{G}_{\beta}(\Gamma_{\alpha}) \right)
\times
\prod_{\beta \leq \gamma < \lambda} 
\left( F_{\gamma}(\Gamma_{\gamma +1}), \mathcal{G}_{\gamma}(\Gamma_{\gamma +1}) \right).
\]
Hence the normaliser tower of $H$ in 
$\Aut^{V^{\mathbb{P}}}(\Gamma)$ terminates in exactly $\beta$
steps.

\end{proof}

\section{Rigid trees} \label{S:rigid}
In this section, we will prove Theorem \ref{T:graph}. Rather than working directly
with graphs, we will find it more convenient to prove the following
analogous theorem for trees. To obtain Theorem \ref{T:graph}, we can then
use one of the standard coding procedures to uniformly convert each tree
$T_{\alpha}$ into a corresponding graph $\Gamma( T_{\alpha} )$. (For
example, we can use the coding of Theorem 5.5.1 \cite{ho}.)

\begin{Thm} \label{T:trees}
It is consistent that for every regular cardinal $\kappa \geq \omega$, there
exists a set $\{ T_{\alpha} \mid \alpha < \kappa^{+} \}$ of 
pairwise nonisomorphic rigid trees of height $\kappa^{+}$ with the
following property. If $E$ is any equivalence relation on
$\kappa^{+}$, then there exists a notion of forcing
$\mathbb{P}$ such that
\begin{enumerate}
\item[(a)] $\mathbb{P}$ preserves cofinalities and cardinalities;
\item[(b)] $\mathbb{P}$ does not adjoin any new $\kappa$-sequences 
of ordinals;
\item[(c)] each tree $T_{\alpha}$ remains rigid in $V^{\mathbb{P}}$;
\item[(d)] $T_{\alpha} \simeq T_{\beta}$ in $V^{\mathbb{P}}$ if{f}
$\alpha \mathrel{E} \beta$.
\end{enumerate}
\end{Thm}

Our proof of Theorem \ref{T:trees} relies heavily on the ideas of Jech \cite{j}.
First we need to introduce some notions from the theory of trees.

\begin{Def} \label{D:tree1}
\begin{enumerate}
\item[(a)] A {\em tree\/} is a partially ordered set 
$\langle T, < \rangle$ such that for every $x \in T$, the set
$\pred_{T}(x) = \{ y \in T \mid y < x \}$ is well-ordered by $<$.
\item[(b)] If $x \in T$, then the {\em height\/} of $x$ in $T$,
denoted $\hit_{T}(x)$, is the order-type of $\pred_{T}(x)$
under $<$.
\item[(c)] If $\alpha$ is an ordinal, then the $\alpha^{th}$
{\em level\/} of $T$ is the set
\[
\Lev_{\alpha}(T) = \{ x \in T \mid \hit_{T}(x) = \alpha \}
\]
and
$T \res \alpha = \underset{\beta < \alpha}{\bigcup} \Lev_{\beta}(T)$.
\item[(d)] A {\em branch\/} of $T$ is a maximal linearly ordered
subset of $T$. The {\em length\/} of a branch $B$ is the order-type
of $B$. An {\em $\alpha$-branch\/} is a branch of length $\alpha$.
\end{enumerate}
\end{Def}

\begin{Def} \label{D:tree2}
Let $\delta$ be an ordinal and let $\lambda$ be a cardinal. A tree $T$
is said to be a {\em $( \delta, \lambda)$-tree\/} if{f}
\begin{enumerate}
\item[(i)] for all $\alpha < \delta$, $0 < \left| \Lev_{\alpha}(T) \right| < \lambda$; and
\item[(ii)] $\Lev_{\delta}(T) = \emptyset$.
\end{enumerate}
A $( \delta, \lambda)$-tree $T$ is {\em normal\/} if each of the following
conditions is satisfied.
\begin{enumerate}
\item[(a)] If $\delta > 0$, then $\left| \Lev_{0}(T) \right| = 1$.
\item[(b)] If $\alpha +1 < \delta$ and $x \in \Lev_{\alpha}(T)$, then there exist
exactly two elements $y_{1}$, $y_{2} \in \Lev_{\alpha +1}(T)$ such that $x < y_{1}$ and
$x < y_{2}$.
\item[(c)] If $\alpha < \beta < \delta$ and $x \in \Lev_{\alpha}(T)$, then there
exists $y \in \Lev_{\beta}(T)$ such that $x < y$.
\item[(d)] Suppose that $\alpha$ is a limit ordinal and $x$, $y \in \Lev_{\alpha}(T)$.
If $\pred_{T}(x) = \pred_{T}(y)$, then $x = y$.
\end{enumerate}
Let $T$ be a $( \delta, \lambda)$-tree. Then the tree $T^{+}$ is an
{\em end-extension\/} of $T$, written $T \lessdot T^{+}$, if
$T^{+} \res \delta = T$. The tree $T^{+}$ is a {\em proper\/}
end-extension if $T \lessdot T^{+}$ and $T \neq T^{+}$.
\end{Def}

\begin{Def} \label{D:tree3}
Let $\kappa \geq \omega$ be a regular cardinal such that
$\kappa^{< \kappa} = \kappa$, and let $\alpha < \kappa^{+}$.
A normal $( \alpha, \kappa^{+})$-tree $T$ is
{\em $<  \kappa$-closed\/} if for each $\beta < \alpha$ such that
$\cf (\beta) < \kappa$ and each increasing sequence of elements
of $T$ 
\[
x_{0} < x_{1} < \dots < x_{\xi} < \dots , \qquad  \xi < \beta ,
\]
such that $x_{\xi} \in \Lev_{\xi}(T)$ for each $\xi < \beta$,
there exists an element $y \in \Lev_{\beta}(T)$ 
such that $\pred_{T}(y) = \{ x_{\xi} \mid \xi < \beta \}$.
\end{Def}

\begin{Lem} \label{L:tree1}
Let $\kappa \geq \omega$ be a regular cardinal such that
$\kappa^{< \kappa} = \kappa$.
\begin{enumerate}
\item[(a)] For each $\alpha < \kappa^{+}$, there exists a
$< \kappa$-closed normal $( \alpha, \kappa^{+})$-tree.
\item[(b)] If $\alpha < \beta < \kappa^{+}$ and $S$ is a
$< \kappa$-closed normal $( \alpha, \kappa^{+})$-tree,
then there exists a
$< \kappa$-closed normal $( \beta, \kappa^{+})$-tree $T$
such that $S \lessdot T$.
\end{enumerate}
\end{Lem}

\begin{proof}
Left to the reader.
\end{proof}

\begin{Lem} \label{L:tree2}
Let $\kappa \geq \omega$ be a regular cardinal such that
$\kappa^{< \kappa} = \kappa$, and let $\alpha < \kappa^{+}$.
If $S$ and $T$ are
$< \kappa$-closed normal $( \alpha, \kappa^{+})$-trees, then
$S \simeq T$. Furthermore, if $\delta +1 \leq \alpha$, then for
each isomorphism 
$\varphi : S \res \delta +1 \to T \res \delta +1$, there exists an
isomorphism $\pi : S \to T$ such that $\varphi \subseteq \pi$.
\end{Lem}

\begin{proof}
If $\alpha < \kappa$, then $S$ and $T$ are both complete binary
trees of height $\alpha$, and so $S \simeq T$. Hence we can suppose
that $\kappa < \alpha < \kappa^{+}$. Thus $|S| = |T| = \kappa$. We
will define an isomorphism 
$\pi = \underset{\xi < \kappa}{\bigcup} \pi_{\xi} : S \to T$ via a
back-and -forth argument.

Suppose that we have defined $\pi_{\xi}$ for some $\xi < \kappa$. 
Assume inductively that there exists a set $\{ B_{i} \mid i \in I \}$
of $\alpha$-branches of $S$ such that
\begin{enumerate}
\item[(i)] $|I| < \kappa$; and
\item[(ii)] $\dom \pi_{\xi} = \underset{i \in I}{\bigcup}B_{i}$.
\end{enumerate}
Let $s$ be any element of $S \smallsetminus \dom \pi_{\xi}$.
Choose an $\alpha$-branch $B$ of $S$ such that $s \in B$. Let
$B = \{ b_{\tau} \mid \tau < \alpha \}$, where 
$b_{\tau} \in B \cap \Lev_{\tau}(S)$. Let $\beta$ be the least ordinal
such that $b_{\beta} \notin \dom \pi_{\xi}$.
First suppose that $\beta = \gamma +1$ is a successor ordinal.
Then there exists an $i \in I$ such that $b_{\gamma} \in B_{i}$.
Let $C_{i} = \pi_{\xi} \left[ B_{i} \right]$. Then there exists a
unique element $c \in \Lev_{\beta}(T) \smallsetminus C_{i}$ such that
$\pi_{\xi}(b_{\gamma}) < c$. Let $C$ be an $\alpha$-branch of $T$
such that $c \in C$; and let $\psi : B \to C$ be the unique
order-preserving bijection. Then 
$\pi_{\xi +1} = \pi_{\xi} \cup \psi$ is a partial isomorphism such that
$s \in \dom \pi_{\xi +1}$. Now suppose that $\beta$ is a limit ordinal.
Since $\{ b_{\tau} \mid \tau < \beta \}$ is covered by the set
$\{ B_{i} \mid i \in I \}$ of branches, it follows that $\cf (\beta) < \kappa$.
Hence there exists an element $c \in \Lev_{\beta}(T)$ such that
$\pred_{T}(c) = \{ \pi_{\xi}(b_{\tau}) \mid \tau < \beta \}$.
Let $C$ be an $\alpha$-branch of $T$
such that $c \in C$; and let $\psi : B \to C$ be the unique
order-preserving bijection. Once again,
$\pi_{\xi +1} = \pi_{\xi} \cup \psi$ is a partial isomorphism such that
$s \in \dom \pi_{\xi +1}$. By a similar argument, if $t$ is any element
of $T \smallsetminus \ran (\pi_{\xi +1})$, then
we can find a partial
isomorphism $\pi_{\xi +2} \supset \pi_{\xi +1}$ such that
$t \in \ran \pi_{\xi +2}$. Hence we can ensure that
$\pi = \underset{\xi < \kappa}{\bigcup} \pi_{\xi}$ is an
isomorphism from $S$ onto $T$.

Finally suppose that $\delta +1 \leq \alpha$ and that
$\varphi : S \res \delta +1 \to T \res \delta +1$ is an isomorphism.
For each $s \in S$ and $t \in T$, define
$S[s] = \{ x \in S \mid s \leq x \}$ and $T[t] = \{ y \in T \mid t \leq y \}$.
Let $\gamma$ be the ordinal such that $\alpha = \delta + \gamma$. Then
for each $s \in \Lev_{\delta}(S)$, both $S[s]$ and $T[ \varphi (t)]$ are 
$< \kappa$-closed normal $( \gamma, \kappa^{+})$-trees; and so
$S[s] \simeq T[ \varphi(s)]$. Hence $\varphi$ can be extended to an
isomorphism $\pi : S \to T$.
\end{proof}

Next we will discuss the notion of forcing
$\mathbb{Q}_{\kappa}$ which adjoins the set 
$\{ T_{\alpha} \mid \alpha < \kappa^{+} \}$ of distinct
pairwise nonisomorphic rigid trees of height $\kappa^{+}$.

\begin{Def} \label{D:poset}
Let $\kappa \geq \omega$ be a regular cardinal such that
$\kappa^{< \kappa} = \kappa$. Then $\mathbb{Q}_{\kappa}$
is the notion of forcing consisting of all conditions
$p = \langle t_{\alpha}^{p} \mid \alpha < \kappa^{+} \rangle$, where
\begin{enumerate}
\item[(a)] each $t_{\alpha}^{p}$ is a $< \kappa$-closed normal
$( \beta_{\alpha}, \kappa^{+})$-tree for some $\beta_{\alpha} < \kappa^{+}$;
and
\item[(b)] there exists an ordinal $\gamma < \kappa^{+}$ such that
$t_{\alpha}^{p} = \emptyset$ for all $\gamma \leq \alpha < \kappa^{+}$.
\end{enumerate}
We define $q \leq p$ if{f} $t^{p}_{\alpha} \lessdot t^{q}_{\alpha}$ for
all $\alpha < \kappa^{+}$.
\end{Def}

Until further notice, we will work with the ground model $M$. Suppose that
$\kappa \geq \omega$ is a regular cardinal such that 
$\kappa^{< \kappa} = \kappa$ and $2^{\kappa} = \kappa^{+}$. Then it is easily 
checked that $\mathbb{Q}_{\kappa}$ is $\kappa$-closed and that
$\left| \mathbb{Q}_{\kappa} \right| = \kappa^{+}$. Hence 
$\mathbb{Q}_{\kappa}$ preserves cofinalities and cardinalities. Let $G$ be an
$M$-generic filter on $\mathbb{Q}_{\kappa}$. For each $\alpha < \kappa^{+}$,
let $T_{\alpha} = \bigcup \{ t_{\alpha}^{p} \mid p \in G \}$.

\begin{Lem} \label{L:rig}
In $M[G]$, $\{ T_{\alpha} \mid \alpha < \kappa^{+} \}$ is a set of distinct
pairwise nonisomorphic rigid trees of height $\kappa^{+}$.
\end{Lem}

\begin{proof}
For each $\alpha < \kappa^{+}$, let $\widetilde{T}_{\alpha}$ be the canonical
$\mathbb{Q}_{\kappa}$-name for $T_{\alpha}$. First suppose that for
some $\alpha < \kappa^{+}$, there exists a nonidentity automorphism $f$
of $T_{\alpha}$ in $M[G]$. Let $\widetilde{f}$ be a 
$\mathbb{Q}_{\kappa}$-name for $f$. Then there exists a condition 
$p \in \mathbb{Q}_{\kappa}$ and an element $a \in t_{\alpha}^{p}$
such that
\[
p \Vdash \widetilde{f} : \widetilde{T}_{\alpha} \to \widetilde{T}_{\alpha}
\text{ is an isomorphism such that } \widetilde{f}(a) \neq a.
\]
Since $\mathbb{Q}_{\kappa}$ is $\kappa$-closed, we can inductively
define a descending sequence of conditions 
$\langle p_{\xi} \mid \xi < \kappa \rangle$ such that
\begin{enumerate}
\item[(1)] $p_{0} = p$;
\item[(2)] $t_{\alpha}^{p_{\xi +1}}$ is a proper end-extension of 
$t_{\alpha}^{p_{\xi}} $; and
\item[(3)] $p_{\xi +1}$ decides $\widetilde{f} \res t_{\alpha}^{p_{\xi}}$.
\end{enumerate}
Let $q$ be the greatest lower bound of the sequence
$\langle p_{\xi} \mid \xi < \kappa \rangle$. Then
$t_{\alpha}^{q}
= \underset{\xi < \kappa}{\bigcup} t_{\alpha}^{p_{\xi}}$,
and so $q$ decides $\widetilde{f} \res t_{\alpha}^{q}$. (In the rest of
this paper, we will refer to the above argument as the 
{\em bootstrap argument\/}.) Note that $t_{\alpha}^{q}$ is
a $< \kappa$-closed normal $( \gamma, \kappa^{+})$-tree
for some ordinal $\gamma$ such that
$\cf (\gamma) = \kappa$. Let $B$ be a $\gamma$-branch
of $t_{\alpha}^{q}$ such that $a \in B$, and let $C$ be the
$\gamma$-branch of $t_{\alpha}^{q}$ such that 
$q \Vdash \widetilde{f} [B] = C$. Then $B \neq C$. Since
$\cf (\gamma) = \kappa$, there exists a
$< \kappa$-closed normal $( \gamma +1, \kappa^{+})$-tree 
$t_{\alpha}^{+}$ such that
\begin{enumerate}
\item[(i)] $t_{\alpha}^{+}$ is a proper end-extension of
$t^{q}_{\alpha}$;
\item[(ii)] there exists $x \in t_{\alpha}^{+}$ such that
$\pred_{ t_{\alpha}^{+}}(x) = B$; and
\item[(iii)] there does not exist $y \in t_{\alpha}^{+}$ such that
$\pred_{ t_{\alpha}^{+}}(y) = C$.
\end{enumerate}
Let $r \leq q$ be a condition such that $t_{\alpha}^{+} 
\lessdot t_{\alpha}^{r}$. Then 
\[
r \Vdash \widetilde{f} \res t_{\alpha}^{q} \text{ cannot be extended to an
automorphism of } t_{\alpha}^{+}.
\]
This is a contradiction.

Now suppose that for some $\alpha < \beta < \kappa^{+}$, there exists an
isomorphism $g : T_{\alpha} \to T_{\beta}$ in $M[G]$. Let $\widetilde{g}$ be
a $\mathbb{Q}_{\kappa}$-name for $g$. Then there exists a condition
$p \in \mathbb{Q}_{\kappa}$ such that
\[
p \Vdash \widetilde{g} : \widetilde{T}_{\alpha} \to \widetilde{T}_{\beta}
\text{ is an isomorphism. }
\]
By the bootstrap argument, there exists a condition $q \leq p$
such that
\begin{enumerate}
\item[(a)] $t_{\alpha}^{q}$ and $t_{\beta}^{q}$ are
$< \kappa$-closed normal $( \gamma, \kappa^{+})$-trees
for some $\gamma$ such that $\cf (\gamma) = \kappa$; 
\item[(b)] $q$ decides $\widetilde{g} \res t_{\alpha}^{q}$.
\end{enumerate}
But then there exist 
$< \kappa$-closed normal $( \gamma +1, \kappa^{+})$-trees 
$t_{\alpha}^{+}$ and
$t_{\beta}^{+}$ such that
\begin{enumerate}
\item[(1)] $t^{+}_{\alpha}$ and $t^{+}_{\beta}$ are proper
end-extensions of $t^{q}_{\alpha}$, $t^{q}_{\beta}$
respectively; and
\item[(2)] $\widetilde{g} \res t_{\alpha}^{q}$ cannot be extended to an 
isomorphism from $t_{\alpha}^{+}$ onto $t_{\beta }^{+}$.
\end{enumerate}
Once again, this yields a contradiction. 
\end{proof}

(A similar argument shows that each $T_{\alpha}$ is a 
$\kappa^{+}$-Suslin tree; cf. the proof of Theorem 48 \cite{j2}.)

Next suppose that $E$ is any equivalence relation on
$\kappa^{+}$. Let $A \subseteq \kappa^{+}$ be the set of
$E$-equivalence class representatives obtained by selecting the
least element of each class.

\begin{Def} \label{D:poset2}
$\mathbb{P}_{E}$ is the notion of forcing in $M[G]$ consisting
of all conditions
$p = \langle f_{\alpha \beta} \mid \alpha < \beta < \kappa^{+} \rangle$ 
such that for some $\gamma < \kappa^{+}$,
\begin{enumerate}
\item[(a)] if $\alpha \in A$, $\beta < \gamma$ and $\alpha \mathrel{E} \beta$,
then there exists $\delta < \kappa^{+}$ such that $f_{\alpha \beta}$ is an
isomorphism from $T_{\alpha} \res \delta +1$ onto $T_{\beta} \res \delta +1$;
\item[(b)] otherwise, $f_{\alpha \beta} = \emptyset$.
\end{enumerate}
The ordering on $\mathbb{P}_{E}$ is the obvious one.
\end{Def}

\begin{Remark} \label{R:why}
Some readers may be wondering why we have introduced the set $A$ of
$E$-equivalence class representatives. Consider the slightly simpler notion
of forcing $\mathbb{P}_{E}^{\prime}$ consisting of all conditions
$p = \langle f_{\alpha \beta} \mid \alpha < \beta < \kappa^{+} \rangle$ 
such that for some $\gamma < \kappa^{+}$,
\begin{enumerate}
\item[(a)] if $\alpha < \beta < \gamma$ and $\alpha \mathrel{E} \beta$,
then there exists $\delta < \kappa^{+}$ such that $f_{\alpha \beta}$ is an
isomorphism from $T_{\alpha} \res \delta +1$ onto $T_{\beta} \res \delta +1$;
\item[(b)] otherwise, $f_{\alpha \beta} = \emptyset$.
\end{enumerate}
Using Lemma \ref{L:tree2}, it is easily seen
that $\mathbb{P}_{E}^{\prime}$ adjoins a generic
isomorphism $g_{\alpha \beta} : T_{\alpha} \to T_{\beta}$
for each $\alpha < \beta < \kappa^{+}$ such that
$\alpha \mathrel{E} \beta$. Fix such a pair $\alpha < \beta$, and
suppose that there exists an ordinal $\gamma$ such that
$\beta < \gamma < \kappa^{+}$ and $\beta \mathrel{E} \gamma$. Then
$g_{\alpha \gamma}$ and $g_{\beta \gamma} \circ g_{\alpha \beta}$
will be {\em distinct\/} isomorphisms from $T_{\alpha}$ onto
$T_{\gamma}$; and so $T_{\alpha}$ will no longer be rigid. The set
$A$ was introduced to deal with precisely this problem. For example,
suppose that $\alpha \in A$. Then the notion of forcing
$\mathbb{P}_{E}$ will only directly adjoin isomorphisms
$g_{\alpha \beta} : T_{\alpha} \to T_{\beta}$ and
$g_{\alpha \gamma} : T_{\alpha} \to T_{\gamma}$. Of course, we can
then obtain an isomorphism from $T_{\beta}$ onto $T_{\gamma}$ by
forming the composition $g_{\alpha \gamma} \circ g_{\alpha \beta}^{-1}$.
\end{Remark}

Let $H$ be an $M[G]$-generic filter on $\mathbb{P}_{E}$. The
following result is an immediate consequence of the discussion in
Remark \ref{R:why}.

\begin{Lem} \label{L:isom}
If $\alpha < \beta < \kappa^{+}$ and $\alpha \mathrel{E} \beta$, then
there exists an isomorphism 
$g_{\alpha \beta} : T_{\alpha} \to T_{\beta}$ in $M[G][H]$.
\end{Lem}
\begin{flushright}
$\square$
\end{flushright}

\begin{Lem} \label{L:point}
$\mathbb{P}_{E}$ preserves cofinalities and cardinalities, and does not
adjoin any new $\kappa$-sequences of ordinals. The following
statements hold in $M[G][H]$.
\begin{enumerate}
\item[(a)] $T_{\alpha}$ is rigid for each $\alpha < \kappa^{+}$.
\item[(b)] If $\alpha < \beta < \kappa^{+}$, then 
$T_{\alpha} \simeq T_{\beta}$ if{f} $\alpha \mathrel{E} \beta$.
\end{enumerate}
\end{Lem}

\begin{proof}
Let $\widetilde{E}$, $\widetilde{A}$ and $\widetilde{\mathbb{P}}_{E}$ be
$\mathbb{Q}_{\kappa}$-names for $E$, $A$ and $\mathbb{P}_{E}$
respectively. Let $\mathbb{R}$ be the subset of 
$\mathbb{Q}_{\kappa} \ast \widetilde{\mathbb{P}}_{E}$ 
consisting of those conditions 
\[
\langle p, \widetilde{q} \rangle =
\langle 
\langle t^{p}_{\alpha} \mid \alpha < \kappa^{+} \rangle,
\langle f_{\alpha \beta} \mid \alpha < \beta < \kappa^{+} \rangle
\rangle
\]
such that for some $\gamma$, $\delta < \kappa^{+}$,
\begin{enumerate}
\item[(1)] $p$ decides $\widetilde{E} \res \gamma \times \gamma$, and hence
$p$ also decides $\widetilde{A} \cap \gamma$; 
\item[(2)] if $\alpha < \gamma$, then $t^{p}_{\alpha}$ is a
$< \kappa$-closed normal $(\delta +1, \kappa^{+})$-tree; and
\item[(3)]
\begin{enumerate}
\item[(i)] if $\alpha < \beta < \gamma$ and 
$p \Vdash \alpha \in \widetilde{A} \text{ and }
\alpha \mathrel{\widetilde{E}} \beta$, then $f_{\alpha \beta}$
is an isomorphism from $t^{p}_{\alpha}$ onto $t^{p}_{\beta}$;
\item[(ii)] otherwise, $f_{\alpha \beta} = \emptyset$.
\end{enumerate}
\end{enumerate}
(Here we are identifying each isomorphism $f_{\alpha \beta}$ with its
canonical $\mathbb{Q}_{\kappa}$-name $\check{f}_{\alpha \beta}$.)

\begin{Claim} \label{C:point}
$\mathbb{R}$ is a dense subset of 
$\mathbb{Q}_{\kappa} \ast \widetilde{\mathbb{P}}_{E}$.
\end{Claim}

\begin{proof}[Proof of Claim \ref{C:point}]
Let $\langle p, \widetilde{q} \rangle =
\langle p,
\langle \widetilde{f}_{\alpha \beta} \mid \alpha < \beta < \kappa^{+} \rangle
\rangle$ be any element of 
$\mathbb{Q}_{\kappa} \ast \widetilde{\mathbb{P}}_{E}$.
Then there exists $p^{\prime} \leq p$ and $\gamma < \kappa^{+}$
such that $p^{\prime}$ forces
\begin{enumerate}
\item[(a)] $ \widetilde{f}_{\alpha \beta} = \emptyset$
for all $\beta \geq \gamma$; and
\item[(b)] if $\alpha \in \widetilde{A}$, $\beta < \gamma$ and 
$\alpha \mathrel{\widetilde{E}} \beta$, then there exists
$\tau < \gamma$ such that $\dom \widetilde{f}_{\alpha \beta}$
is a $< \kappa$-closed normal $(\tau +1, \kappa^{+})$-tree.
\end{enumerate}
Since $\mathbb{Q}_{\kappa}$ is $\kappa$-closed, there exists
$r \leq p^{\prime}$ such that
\begin{enumerate}
\item[(c)] $r$ decides $\widetilde{E} \res \gamma \times \gamma$,
and hence $r$ also decides $\widetilde{A} \cap \gamma$;
\item[(d)] there exists $\delta \geq \gamma$ such that 
$t^{r}_{\alpha}$ is a $< \kappa$-closed normal
$(\delta +1, \kappa^{+})$-tree for each $\alpha < \gamma$; and
\item[(e)] if $\alpha < \beta < \gamma$ and 
$r \Vdash \alpha \in \widetilde{A} \text{ and }
\alpha \mathrel{\widetilde{E}} \beta$, then there exists
$\tau < \gamma$ and an isomorphism 
$f_{\alpha \beta} : t^{r}_{\alpha} \res \tau +1 \to t^{r}_{\beta} \res \tau +1$
such that $r \Vdash \widetilde{f}_{\alpha \beta} = f_{\alpha \beta}$.
\end{enumerate}
By Lemma \ref{L:tree2}, if $\alpha < \beta < \gamma$ and 
$r \Vdash \alpha \in \widetilde{A} \text{ and }
\alpha \mathrel{\widetilde{E}} \beta$, then there exists an
isomorphism $g_{\alpha \beta} : t^{r}_{\alpha} \to t^{r}_{\beta}$ such
that $f_{\alpha \beta} \subset g_{\alpha \beta}$. Let
$ g_{\alpha \beta} = \emptyset$ for all other pairs
$\alpha < \beta < \kappa^{+}$. Then
$\langle r, 
\langle g_{\alpha \beta} \mid \alpha < \beta < \kappa^{+} \rangle
\rangle \in \mathbb{R}$ is a strengthening of 
$\langle p, \widetilde{q} \rangle$.
\end{proof}

Thus the forcing notions 
$\mathbb{Q}_{\kappa} \ast \widetilde{\mathbb{P}}_{E}$
and $\mathbb{R}$ are equivalent. It is easily checked that $\mathbb{R}$ is
$\kappa$-closed and that $\left| \mathbb{R} \right| = \kappa^{+}$. Hence
$\mathbb{R}$ preserves cofinalities and cardinalities, and does not
adjoin any new $\kappa$-sequences of ordinals. It follows that the same
is true of $\mathbb{P}_{E} \in M[G]$.

Now suppose that for some $\mu < \kappa^{+}$, there exists a nonidentity
automorphism $\varphi$ of $T_{\mu}$ in $M[G][H]$. Let
$\widetilde{\varphi}$ be an $\mathbb{R}$-name for $\varphi$. Then there
exists a condition $\langle p, \widetilde{q} \rangle \in \mathbb{R}$ and an
element $a \in t^{p}_{\mu}$ such that
\[
\langle p, \widetilde{q} \rangle \Vdash
\widetilde{\varphi} : \widetilde{T}_{\mu} \to \widetilde{T}_{\mu}
\text{ is an automorphism such that } 
\widetilde{\varphi}(a) \neq a.
\]
Since $\mathbb{R}$ is $\kappa$-closed, we can inductively define
a descending sequence of conditions
$\langle p_{\xi}, \widetilde{q}_{\xi} \mid \xi < \kappa \rangle$
such that
\begin{enumerate}
\item[(i)] $\langle p_{0} , \widetilde{q}_{0} \rangle =
\langle p , \widetilde{q} \rangle$;
\item[(ii)] $t_{\mu}^{p_{\xi +1}}$ is a proper end-extension of
$t_{\mu}^{p_{\xi}}$; and
\item[(iii)] $\langle p_{\xi +1}, \widetilde{q}_{\xi +1} \rangle$
decides $\widetilde{\varphi} \res t_{\mu}^{p_{\xi}}$.
\end{enumerate}
Let $t_{\mu} = \underset{\xi < \kappa}{\bigcup}t_{\mu}^{p_{\xi}}$;
and let $\psi : t_{\mu} \to t_{\mu}$ be the nonidentity 
automorphism such that for each $\xi < \kappa$,
$\langle p_{\xi +1}, \widetilde{q}_{\xi +1} \rangle
\Vdash \widetilde{\varphi} \res t_{\mu}^{p_{\xi}}
\subset \psi$. Note that $t_{\mu}$ is a
$< \kappa$-closed normal $(\eta , \kappa^{+})$-tree
for some $\eta$ such that $\cf (\eta) = \kappa$. Arguing as in the
proof of Lemma \ref{L:rig}, we see that there exists a
$< \kappa$-closed normal $(\eta +1, \kappa^{+})$-tree 
$t^{+}_{\mu} \supset t_{\mu}$ such that $\psi$ cannot
be extended to an automorphism of $t^{+}_{\mu}$. 
But then the following claim yields a contradiction.

\begin{Claim} \label{C:detail}
There exists a condition
$\langle p_{\kappa}, \widetilde{q}_{\kappa} \rangle
\in \mathbb{R}$ such that
\begin{enumerate}
\item[(1)] $\langle p_{\xi}, \widetilde{q}_{\xi} \rangle \leq
\langle p_{\kappa}, \widetilde{q}_{\kappa} \rangle$ for all
$\xi < \kappa$; and
\item[(2)] $t_{\mu}^{+} \lessdot t_{\mu}^{p_{\kappa}}$.
\end{enumerate}
\end{Claim}

\begin{proof}[Proof of Claim \ref{C:detail}]
For each $\alpha < \kappa^{+}$, let
$t_{\alpha} = \underset{\xi < \kappa}{\bigcup}t_{\alpha}^{p_{\xi}}$.
In order to construct a suitable condition
$\langle p_{\kappa}, \widetilde{q}_{\kappa} \rangle
\in \mathbb{R}$, we must be able to simultaneously solve the 
following extension problems. For various pairs of ordinals
$\alpha < \beta < \kappa^{+}$, we are given an isomorphism
$f_{\alpha \beta} : t_{\alpha} \to t_{\beta}$; and we must find
suitable extensions $t_{\alpha}^{p_{\kappa}}$ and
$t_{\beta}^{p_{\kappa}}$ of $t_{\alpha}$, $t_{\beta}$
such that $f_{\alpha \beta}$ extends to an isomorphism of
$t_{\alpha}^{p_{\kappa}}$ onto $t_{\beta}^{p_{\kappa}}$.
Of course, the most difficult cases are when either
$\alpha = \mu$ or $\beta = \mu$; for then we have the additional
requirement that $t_{\mu}^{+} \lessdot t_{\mu}^{p_{\kappa}}$.
However, for each such pair of ordinals 
$\alpha < \beta < \kappa^{+}$, there exists $\xi < \kappa$ such
that $p_{\xi} \Vdash \alpha \in \widetilde{A}$. Consequently
$f_{\alpha \beta}$ is the only isomorphism which needs to be
considered when extending $t_{\beta}$ to
$t_{\beta}^{p_{\kappa}}$; and so there are no conflicts.
\end{proof}

A similar argument shows that if $\alpha < \beta < \kappa^{+}$
and $\alpha$, $\beta$ are not $E$-equivalent, then $T_{\alpha}$
and $T_{\beta}$ remain nonisomorphic in $M[G][H]$.
\end{proof}

Finally we will use a reverse Easton iteration to complete the proof 
of Theorem \ref{T:trees}.  (Clear
accounts of reverse Easton forcing can be found in Baumgartner
\cite{b} and Menas \cite{m}.) 
Let $V_{0}$ be a transitive model of $ZFC + GCH$. We define a
sequence of posets 
$\langle \mathbb{P}_{\alpha} \mid \alpha \in On \rangle$ inductively
as
follows. \\
{\bf Case 1.} If $\alpha = 0$, then $\mathbb{P}_{0}$ is the trivial
poset such that $\left| \mathbb{P}_{0} \right| = 1$. \\
{\bf Case 2.} If $\alpha$ is a limit ordinal which is not inaccessible,
then $\mathbb{P}_{\alpha}$ is the inverse limit of
$\langle \mathbb{P}_{\beta} \mid \beta < \alpha \rangle$. \\
{\bf Case 3.} If $\alpha$ is inaccessible, then
$\mathbb{P}_{\alpha}$ is the direct limit of
$\langle \mathbb{P}_{\beta} \mid \beta < \alpha \rangle$. \\
{\bf Case 4.} Finally suppose that $\alpha = \gamma +1$. If
$\gamma = \kappa \geq \omega$ is a regular cardinal, then
$\mathbb{P}_{\alpha} = 
\mathbb{P}_{\kappa} \ast \widetilde{\mathbb{Q}}_{\kappa}$,
where $\widetilde{\mathbb{Q}}_{\kappa} \in
V_{0}^{\mathbb{P}_{\kappa}}$ is the notion of forcing introduced
in Definition \ref{D:poset}. Otherwise,
$\mathbb{P}_{\alpha} = 
\mathbb{P}_{\gamma} \ast \mathbb{P}_{0}$.

Let $\mathbb{P}_{\infty}$ be the direct limit of
$\langle \mathbb{P}_{\alpha} \mid \alpha \in On \rangle$.
For each $\alpha \in On$, let $\mathbb{P}_{\alpha \infty}$ be the
canonically chosen class in $V_{0}^{\mathbb{P}_{\alpha}}$
such that $\mathbb{P}_{\infty} \simeq \mathbb{P}_{\alpha}
\ast \mathbb{P}_{\alpha \infty}$. Let the class
$G \subseteq \mathbb{P}_{\infty}$ be $V_{0}$-generic; and for
each $\alpha \in On$, let $G_{\alpha} = G \cap \mathbb{P}_{\alpha}$.
Let $V = V_{0}[G] = \underset{\alpha \in On}{\bigcup}V_{0}[G_{\alpha}]$.
Then a routine argument yields the following result.

\begin{Lem} \label{L:easton}
\begin{enumerate}
\item[(a)] $\mathbb{P}_{\infty}$ preserves cofinalities and cardinalities.
\item[(b)] $V$ is a model of $ZFC + GCH$.
\end{enumerate}
\end{Lem}
\begin{flushright}
$\square$
\end{flushright}

Let $\kappa \geq \omega$ be any regular cardinal; and let
$\{ T_{\alpha} \mid \alpha < \kappa^{+} \} \in
V_{0}[G_{\kappa +1}]$ 
be the set of trees
which is adjoined by $\mathbb{Q}_{\kappa}$ at the
$\kappa^{th}$ stage of the iteration. Since 
$\mathbb{P}_{\kappa +1 \infty}$ is $\kappa^{+}$-closed,
it follows that $\{ T_{\alpha} \mid \alpha < \kappa^{+} \}$ 
remains a set of pairwise nonisomorphic
rigid trees in $V$. Now let $E \in V$ be any equivalence relation
on $\kappa^{+}$, and let $\mathbb{P}_{E}$ be the corresponding
notion of forcing, which was introduced in Definition \ref{D:poset2}.
Again using the fact that
$\mathbb{P}_{\kappa +1 \infty}$ is $\kappa^{+}$-closed, we see that
$E$, $\mathbb{P}_{E} \in V_{0}[G_{\kappa +1}]$. We have already
shown that $\mathbb{P}_{E}$ has the appropriate properties in
$ V_{0}[G_{\kappa +1}]$. Thus it only remains to prove that these
properties are preserved in $V$.

\begin{Lem} \label{L:preserve}
In $V$, $\mathbb{P}_{E}$ preserves cofinalities and cardinalities, 
and does not adjoin any new $\kappa$-sequences of ordinals. 
The following statements hold in $V^{\mathbb{P}_{E}}$. 
\begin{enumerate}
\item[(a)] $T_{\alpha}$ is rigid for each $\alpha < \kappa^{+}$.
\item[(b)] If $\alpha < \beta < \kappa^{+}$, then 
$T_{\alpha} \simeq T_{\beta}$ if{f} $\alpha \mathrel{E} \beta$.
\end{enumerate}
\end{Lem}

\begin{proof}
Since $\left| \mathbb{P}_{E} \right| = \kappa^{+}$, $\mathbb{P}_{E}$
preserves cofinalities and cardinalities greater than $\kappa^{+}$. The
remaining parts of the lemma correspond to combinatorial properties
of $\mathbb{P}_{E}$ which are preserved under
$\kappa^{+}$-closed forcing. For example,
suppose that $p \in \mathbb{P}_{E}$ satisfies
\[
p \Vdash \widetilde{f} : T_{\alpha} \to T_{\alpha}
\text{ is an automorphism. }
\]
We can assume that $\widetilde{f}$ is a nice $\mathbb{P}_{E}$-name;
ie. that 
$\widetilde{f} = \bigcup \{ \{ \langle s,t \rangle \} \times A_{s,t} \mid
s, t \in T_{\alpha} \}$,
where each $A_{s,t}$ is an antichain of $\mathbb{P}_{E}$. Then
$\widetilde{f} \in V_{0}[G_{\kappa +1}]$, and so there exists
$q \leq p$ such that
\[
q \Vdash \widetilde{f}(t) = t \text{ for all } t \in T_{\alpha}.
\]
Hence $T_{\alpha}$ is rigid in $V^{\mathbb{P}_{E}}$.
\end{proof}

This completes the proof of Theorem \ref{T:trees}.

\end{document}